\documentclass[12pt]{article}
\font\smallit=cmti10
\font\smalltt=cmtt10
\font\smallrm=cmr9

 \usepackage{amsmath}
 \usepackage{amsthm}
 \usepackage{amssymb,latexsym}

 at 10 true pt

\def\be#1{ \begin{equation}\label{#1} }

\def\bas{\begin{align*}}
\def\eas{\end{align*}}
\def\bi{\begin{itemize}}
\def\ei{\end{itemize}}

\def\Z{{\hbox{\bf Z}}}

\def\emph#1{{\it #1}}
\def\textbf#1{{\bf #1}}

\def\ep{{\epsilon}}

\def \ep{\varepsilon}

\theoremstyle{plain}
  \newtheorem{theorem}[subsection]{Theorem}
  \newtheorem{question}[subsection]{Question}

  \newtheorem{lemma}[subsection]{Lemma}

  \newtheorem{claim}[subsection]{Claim}
\theoremstyle{remark}

\theoremstyle{definition}

\makeatletter 

\renewcommand\section{\@startsection {section}{1}{\z@}%
                                   {-30pt  \@plus -1ex \@minus -.2ex}%
                                   {2.3ex \@plus.2ex}%
                                   {\normalfont\normalsize\bfseries}}

\renewcommand\subsection{\@startsection{subsection}{2}{\z@}%
                                     {-3.25ex\@plus -1ex \@minus -.2ex}%
                                     {1.5ex \@plus .2ex}%
                                     {\normalfont\normalsize\bfseries}}

\renewcommand{\@seccntformat}[1]{\csname the#1\endcsname. } 

\makeatother

\begin{document}
\vspace*{-40pt}
\centerline{\smalltt INTEGERS: \smallrm
ELECTRONIC JOURNAL OF COMBINATORIAL NUMBER THEORY \smalltt 9
(2009), \#A02}
\vskip 40pt

\begin{center}
\uppercase{\bf On two-point configurations in a random set}
\vskip 20pt
{\bf Hoi H. Nguyen\footnote{This work was written while the author was supported by a DIMACS summer research fellowship, 2008.}}\\
{\smallit Department of Mathematics, Rutgers University, Piscataway, NJ 08854, USA}\\
{\tt hoi@math.rutgers.edu}\\
\end{center}
\vskip 30pt
\centerline{\smallit Received: 8/2/08,
Accepted: 1/3/09, Published: 1/8/09} 
\vskip 30pt

\centerline{\bf Abstract}

We show that with high probability a random subset of
$\{1,\dots,n\}$ of size $\Theta(n^{1-1/k})$ contains two elements
$a$ and $a+d^k$, where $d$ is a positive integer. As a consequence,
we prove an analogue of the S\'ark\"ozy-F\"urstenberg theorem for a
random subset of $\{1,\dots,n\}$.

\pagestyle{myheadings}
\markright{\smalltt INTEGERS: \smallrm ELECTRONIC JOURNAL OF COMBINATORIAL NUMBER THEORY \smalltt 9 (2009), \#A03\hfill} 

\thispagestyle{empty}
\baselineskip=15pt
\vskip 30pt

\section{Introduction}

\noindent Let $\wp$ be a general additive configuration,
$\wp=(a,a+P_1(d),\dots,a+P_{k-1}(d))$, where $P_i\in \Z[d]$ and
$P_i(0)=0$. Let $[n]$ denote the set of positive integers up to $n$. A natural question is:

\begin{question}\label{question:1}
How is $\wp$ distributed in $[n]$?
\end{question}

Roth's theorem \cite{Ro} says that for $\delta>0$ and sufficiently
large $n$, any subset of $[n]$ of size $\delta n$ contains a
nontrivial instance of $\wp = (a,a+d,a+2d)$ (here nontrivial means
$d\neq 0$). In 1975, Szemer\'edi \cite{Sz} extended Roth's theorem
for general linear configurations $\wp=(a,a+d,\dots, a+(k-1)d)$. For
a configuration of type $\wp=(a,a+P(d))$, S\'ark\"ozy \cite{Sa} and
F\"urstenberg \cite{Fu} independently discovered a similar
phenomenon.

\begin{theorem}[S\'ark\"ozy-F\"urstenberg theorem, quantitative
version]\em{\cite[Theorem 3.2]{Ta-Zi},\cite[Theorem
3.1]{Ha-La}}\label{theorem:Sa-Fu} Let $\delta$ be a fixed positive
real number, and let $P$ be a polynomial of integer coefficients
satisfying $P(0)=0$. Then there exists an integer $n=n(\delta,P)$
and a positive constant $c(\delta,P)$ with the following property.
If $n\ge n(\delta,P)$ and $A\subset [n]$ is any subset of
cardinality at least $\delta n$, then \vspace*{-15pt}
\begin{itemize}

\item $A$ contains a nontrivial instance of
$\wp$.

\item $A$ contains at least $c(\delta,P)|A|^2n^{1/\deg(P)-1}$ instances of $\wp=(a,a+P(d))$.
\end{itemize}
\end{theorem}

 In 1996, Bergelson and Leibman \cite{Be-Le} extended this
result for all configurations $\wp=(a,a+P_1(d),\dots, P_{k-1}(d))$, where $P_i \in \Z[d]$ and $P_i(0)=0$ for all $i$.

 Following Question \ref{question:1}, one may
consider the distribution of $\wp$ in a ``pseudo-random" set.

\begin{question}\label{question:2}
Does the set of primes contain a nontrivial instance of $\wp$? How
is $\wp$ distributed in this set?
\end{question}

The famous Green-Tao theorem \cite{Gr-Ta} says that any subset of
positive upper density of the set of primes contains a nontrivial
instance of $\wp=(a,a+d,\dots, a+(k-1)d)$ for any $k$. This
phenomenon also holds for more general configurations
$(a,a+P_1(d),\dots, a+P_{k-1}(d))$, where $P_i\in \Z[d]$ and
$P_i(0)=0$ for all $i$ (cf. \cite{Ta-Zi}).

 The main goal of this note is to consider a similar
question.

\begin{question}\label{question:3}
How is $\wp$ distributed in a typical random subset of $[n]$?
\end{question}

Let $\wp$ be an additive configuration and let $\delta$ be a fixed
positive real number. We say that a set $A$ is $(\delta,\wp)$-dense
if any subset of cardinality at least $\delta |A|$ of $A$ contains a
nontrivial instance of $\wp$. In 1991, Kohayakawa-{\L}uczak-R\"odl
\cite{Ko-Lu-Ro} showed the following result.

\begin{theorem}\label{theorem:Ko-Lu-Ro}
Almost every subset $R$ of $[n]$ of cardinality $|R|=r\gg_\delta
n^{1/2}$ is $(\delta,(a,a+d,a+2d))$-dense.
\end{theorem}

 The assumption $r\gg_\delta n^{1/2}$ is tight, up to a constant factor. Indeed, a typical random subset $R$ of $[n]$ of cardinality $r$ contains about $\Theta(r^3/n)$ three-term arithmetic progressions. Hence, if $(1-\delta) r \gg r^3/n$, then there is a subset of $R$ of cardinality $\delta r$ which does not contain any nontrivial 3-term arithmetic progression.

 Motivated by Theorem \ref{theorem:Ko-Lu-Ro}, {\L}aba and Hamel \cite{Ha-La} studied the distribution of $\wp=(a,a+d^k)$ in a typical random subset of $[n]$, as follows.

\begin{theorem}\label{theorem:Ha-La} Let $k\ge 2$ be an integer. Then there exists a positive real number $\ep(k)$ with the
following property. Let $\delta$ be a fixed positive real number,
then almost every subset
  $R$ of $[n]$ of cardinality $|R|=r\gg_\delta n^{1-\ep(k)}$ is $(\delta,(a,a+d^k))$-dense.
\end{theorem}

 It was shown that $\ep(2)=1/110$, and $\ep(3)\gg
\ep(2)$, etc. Although the method used in \cite{Ha-La} is strong, it
seems to fall short of obtaining relatively good estimates for
$\ep(k)$. On the other hand, one can show that $\ep(k) \le 1/k$.
Indeed, a typical random subset of $[n]$ of size $r$ contains
$\Theta(n^{1+1/k}r^2/n^2)$ instances of $(a,a+d^k)$. Thus if
$(1-\delta) r \gg n^{1+1/k}r^2/n^2 $ (which implies $r\ll_{\delta}
n^{1-1/k}$) then there is a subset of size $\delta r$ of $R$ which
does not contain any nontrivial instance of $(a,a+d^k)$.

 In this note we shall sharpen Theorem \ref{theorem:Ha-La} by showing that $\ep(k)=1/k$.

\begin{theorem}[Main theorem]\label{theorem:integer}
Almost every subset $R$ of $[n]$ of size $|R|=r\gg_\delta n^{1-1/k}$
is  $(\delta,(a,a+d^k))$-dense.
\end{theorem}

 Our method to prove Theorem \ref{theorem:integer} is elementary. We will invoke a combinatorial lemma and the quantitative S\'ark\"ozy-F\"urstenberg theorem (Theorem \ref{theorem:Sa-Fu}). As
the reader will see later on, the method also works for more general
configurations $(a,a+P(d))$, where $P\in \Z[d]$ and $P(0)=0$.

\section{A Combinatorial Lemma}\label{Section:graph}
\vspace*{-10pt}
 Let $G(X,Y)$ be a bipartite graph. We denote the
number of edges going through $X$ and $Y$ by $e(X,Y)$. The average
degree $\bar{d}(G)$ of $G$ is defined to be $e(X,Y)/(|X||Y|)$.

\begin{lemma}\label{lemma:graph}
Let $\{G=G([n],[n])\}_{n=1}^\infty$ be a sequence of bipartite
graphs. Assume that for any $\ep>0$ there exist an integer $n(\ep)$
and a number $c(\ep)>0$ such that $e(A,A)\ge c(\ep) |A|^2
\bar{d}(G)/n$ for all $n\ge n(\ep)$ and all $A\subset [n]$
satisfying $|A|\ge \ep n$. Then for any $\alpha >0$ there exist an
integer $n(\alpha)$ and a number $C(\alpha)>0$ with the following
property. If one chooses a random subset $S$ of $[n]$ of cardinality
$s$, then the probability of $G(S,S)$ being empty is at most
$\alpha^{s}$, providing that $|S| = s \ge C(\alpha)n/\bar{d}(G)$ and
$n\ge n(\alpha)$.
\end{lemma}

\begin{proof} For short we denote the ground set $[n]$ by $V$. We shall view $S$ as an ordered random subset, whose elements will be chosen in order, $v_1$ first and $v_s$ last. We shall verify the lemma within  this probabilistic model. Deduction of the original model follows easily.

 For $1\le k\le s-1$, let $N_k$ be the set of neighbors of the first $k$ chosen vertices, i.e., $N_k=\{v\in V, (v_i,v)\in E(G) \mbox{ for some } i\le k \}$. Since $G(S,S)$ is empty, we have $v_{k+1}\notin N_k$. Next, let $B_{k+1}$ be the set of possible choices
for $v_{k+1}$ (from $V \backslash \{v_1,\dots,v_k\}$) such that $N_{k+1}\backslash N_k\le
c(\ep)\ep \bar{d}(G)$, where $\ep$ will be chosen to be small enough ($\ep =\alpha^2/6$ is fine) and $c(\ep)$ is the constant from Lemma \ref{lemma:graph}. We observe the following.

\begin{claim}\label{Claim:degree2:graph:A,A}
$|B_{k+1}|\le \ep |V|$.
\end{claim}
\vspace*{-10pt}
To prove this claim, we assume for contradiction that $|B_{k+1}|\ge \ep |V|=\ep n$. Since $B_{k+1}\cap
  N_k=\emptyset$, we have
$e(B_{k+1},B_{k+1}) \le e(B_{k+1},V\backslash N_k) \le c(\ep) \ep
  \bar{d}(G)|B_{k+1}| < c(\ep) |B_{k+1}|^2\bar{d}(G)/n$.
This contradicts the property of $G$ assumed in Lemma
\ref{lemma:graph}, provided that $n$ is large enough.

 Thus we conclude that if $G(S,S)$ is empty then
$|B_{k+1}|\le \ep |V|$ for $1\le k\le s-1$.

 Now let $s$ be sufficiently large, say
$s\ge 2(c(\ep)\ep)^{-1} n/\bar{d}(G)$, and assume that the vertices
$v_1,\dots,v_s$  have been chosen. Let $s'$ be the number of
vertices $v_{k+1}$ that do not belong to $B_{k+1}$. Then we have

$$n\ge |N_s|\ge  \sum_{v_{k+1}\notin B_{k+1}}|N_{k+1}\backslash N_k| \ge s' c(\ep) \ep \bar{d}(G).$$

Hence, $s'\le (c(\ep) \ep)^{-1}n/\bar{d}(G) \le s/2.$

 As a result, there are $s-s'$ vertices $v_{k+1}$ that belong to
$B_{k+1}$. But since $|B_{k+1}|\le \ep n$, we see that the number of
subsets $S$ of $V$ such that $G(S,S)$ is empty is bounded by
$$\sum_{s'\le s/2}\binom{s}{s'}n^{s'}(\ep n)^{s-s'} \le (6\ep)^{s/2}n(n-1)\dots (n-s+1)\le \alpha^s n(n-1)\dots (n-s+1),$$
thereby completing the proof.
\end{proof}

\section{Proof of Theorem \ref{theorem:integer}}

First, we define a bipartite graph $G$ on $[n]\times [n] = V_1 \times V_2$ by connecting
$u\in V_1$  to $v \in V_2$ if $v-u=d^k$
for some integer $d\in [1,n^{1/k}]$. Notice that $\bar{d}(G) \approx Cn^{1/k}$ for some absolute
constant $C$.

 Let us restate the S\'ark\"ozy-F\"urstenberg theorem
(Theorem \ref{theorem:Sa-Fu}, for $P(d)=d^k$) in terms of the graph
$G$.

\begin{theorem}\label{theorem:S-F}
Let $\ep>0$ be a positive constant. Then there exists a positive
integer $n(\ep,k)$ and a positive constant $c(\ep,k)$ such that
$e(A,A)\ge c(\ep,k) |A|^2 n^{1/k-1}$ for all $n\ge n(\ep,k)$ and all
$A\subset [n]$ satisfying $|A|\ge \ep n$.
\end{theorem}

Now let $S$ be a subset of $[n]$ of size  $s$. We call $S$ {\it bad}
if it does not contain any nontrivial instance of $(a,a+d^k)$. In
other words, $S$ is bad if $G(S,S)$ contains  no edges.
 By Lemma \ref{lemma:graph} and Theorem \ref{theorem:S-F}, the number of bad
subsets of $[n]$ is at most $\alpha^s \binom{n}{s}$, provided that $s \ge
C(\alpha)n/\bar{d}(G)$. This condition is satisfied if we assume that

$$s\ge 2C(\alpha)C^{-1}
n^{1-1/k}.$$

Next, let $r=s/\delta$ and consider a random subset $R$ of $[n]$ of
size $r$. The probability that $R$ contains a bad subset of size $s$
is at most
$$\alpha^s \binom{n}{s}\binom{n-s}{r-s}/\binom{n}{r}=o(1),$$
provided that $\alpha = \alpha(\delta)$ is small enough.

 To finish the proof, we note that if $R$ does not contain any bad
subset of size $\delta r$, then $R$ is $(\delta,(a,a+d^k))$-dense.


\begin{thebibliography}{9} \footnotesize

\bibitem{Be-Le} {V. Bergelson} and {A. Leibman}, {\it Polynomial extensions of Van Der Waerden's and Szemer\'edi's theorems}, J. Amer. Math. Soc. 9(1996), no. 3, 725-753.

\bibitem{Fu}{H. F\"urstenberg} {\it Recurrence in ergodic theory
  and combinatorial number theory}, Princeton University Press,
  Princeton (1981).


\bibitem{Gr-Ta} {B. J. Green} and {T. Tao}, {\it Primes contain arbitrarily long arithmetic progression}, to appear in Ann. Math.


\bibitem{Ha-La} {M. Hamel} and {I. {\L}aba}, {\it Arithmetic structures in random sets}, Electronic Journal of Combinatorial Number Theory 8 (2008).

\bibitem{Ko-Lu-Ro}{Y. Kohayakawa, T. {\L}uczak}, and {V. R\"odl}, {\it Arithmetic progressions of length three in subsets on a random sets}, Acta Arith. 75(1996), 133-163.

\bibitem{Ro} {K. F. Roth}, {\it On certain sets of integers}, J. London Math Soc. 28(1953), 245-252.

\bibitem{Sa} {A. S\'ark\"ozy}, {\it On difference sets of integres
III}, Acta Math. Sci. Hungar., 31, (1978), 125-149.

\bibitem{Sz} {E. Szemer\'edi}, {\it On set of integers containing no $k$ elements in arithmetic progressions}, Acta Arith. 27 (1975), 299-345.

\bibitem{Ta-Zi} {T. Tao} and {T. Ziegler} {\it The primes contain arbitrarily long polinomial progressions}, to appear in Acta Math.


\end{thebibliography}
\end{document}